\documentclass[11pt]{article}
\usepackage{amsmath,amsthm,amsfonts,latexsym,amscd}

\swapnumbers
\theoremstyle{plain}
\newtheorem{theorem}{Theorem}[section]
\newtheorem{lemma}[theorem]{Lemma}
\newtheorem{corollary}[theorem]{Corollary}
\newtheorem{proposition}[theorem]{Proposition}
\newtheorem{conjecture}[theorem]{Conjecture}

\theoremstyle{definition}
\newtheorem{definition}[theorem]{Definition}
\newtheorem{example}[theorem]{Example}

\newtheorem{situation}[theorem]{Situation}

\theoremstyle{remark}
\newtheorem{remark}[theorem]{Remark}

\usepackage{amsmath,amsfonts,latexsym}

\newcommand{\R}{\mathbb{R}}
\newcommand{\reals}{\mathbb{R}}
\newcommand{\C}{\mathbb{C}}
\newcommand{\complexs}{\mathbb{C}}
\newcommand{\N}{\mathbb{N}}
\newcommand{\naturals}{\mathbb{N}}
\newcommand{\Z}{\mathbb{Z}}
\newcommand{\integers}{\mathbb{Z}}

\newcommand{\1}{\mathbf{1}}  % if this does not work: \boldsymbol{1}
\DeclareMathOperator{\id}{id}

\newcommand{\boundedops}{\mathcal{B}}

\newcommand{\abs}[1]{\left\lvert#1\right\rvert} %absolute value
\newcommand{\norm}[1]{\left\lVert#1\right\rVert}
\newcommand{\generate}[1]{\langle#1\rangle}
\newcommand{\tensor}{\otimes}

\newcommand{\disjointunion}{\amalg}
\newcommand{\subgroup}{<}

\DeclareMathOperator{\tr}{tr}
\DeclareMathOperator{\dirlim}{dirlim}
\DeclareMathOperator{\invlim}{invlim}

\newcommand{\forget}[1]{}
\newcommand{\innerprod}[1]{\langle #1 \rangle}

{\catcode`@=11\global\let\c@equation=\c@theorem}

\allowdisplaybreaks[1]
\newcommand{\RAgroups}{\mathcal G}
\newcommand{\NeumannN}{\mathcal N}

\DeclareMathOperator{\diag}{diag}

\DeclareMathOperator{\Det}{det}

\begin{document}
%%%%%%%%%%%%%%%%%%%%%%% immer wieder dran denken %%%%%%%%%%%%%%%%%%%
\date{last compiled: \today, last edited after: August 10, 2011}
%%%%%%%%%%%%%%%%%%%%%%%%%%%%%%%%%%%%%%%%%%%%%%%%%%%%%%%%%%%%%%%%%%%%

\title{$L^2$-determinant class and approximation of $L^2$-Betti numbers}
\forget{\author{Thomas Schick}
\address{Mathematisches Institut\\
Georg-August Universit\"at G\"ottingen\\
Bunsenstr.~3\\
         37073 G\"ottingen }
\email{schick@uni-math.gwdg.de}
\urladdr{http://www.uni-math.gwdg.de/schick/}
\keywords{$L^2$-determinant, $L^2$-Betti numbers, approximation,
  $L^2$-torsion,  homotopy invariance}
\subjclass{Primary: 58G26, Secondary: 55N25, 55P29}
}
\author{Thomas Schick\thanks{
e-mail: schick@uni-math.gwdg.de\protect\\
www:~http://www.uni-math.gwdg.de/schick/
}\\Mathematisches Institut ---
Georg-August Universit\"at G\"ottingen\\
Bunsenstr. 3 ---
         37073 G\"ottingen }
\maketitle

\begin{abstract}
  A standing conjecture in $L^2$-cohomology is that every finite
  $CW$-complex $X$ is of $L^2$-determinant class. In this
  paper, we prove this  whenever the fundamental group belongs to a
  large class $\RAgroups$ of groups containing e.g.~all extensions of
  residually finite groups with amenable quotients, all 
  residually amenable groups and free products of these. If, in
  addition, $X$ is $L^2$-acyclic, we also
  prove that the $L^2$-determinant is a homotopy invariant.
  Even in
  the known cases, our proof of homotopy invariance is much shorter and easier than the
  previous ones.
  Under suitable conditions  we
  give new approximation formulas for $L^2$-Betti numbers. 
\forget{
  The paper generalizes results of L\"uck;
  Dodziuk, Mathai, Rothenberg; and Clair.
}

Errata are added, rectifying some unproved statements about ``amenable
extension'': throughout, amenable extensions should be extensions with
\emph{normal} subgroups.

Keywords: {$L^2$-determinant, $L^2$-Betti numbers, approximation,
  $L^2$-torsion,  homotopy invariance}\\
MSC: 58G26 (Primary),  55N25, 55P29 (Secondary)

\end{abstract}

%\maketitle

\section{Introduction}

For a finite $CW$-complex $X$ with fundamental group $\pi$,
$L^2$-invariants of the universal covering $\tilde X$ are defined in
terms of the combinatorial Laplacians $\Delta_*$ on $C^*_{(2)}(\tilde
X)=C^*_{cell}(\tilde X)\tensor_{\integers\pi} l^2(\pi)$, which, after
the choice of a cellular base, is a finite direct sum of copies of
$l^2(\pi)$. $\Delta_p= (c_p\tensor\id)^*(c_p\tensor\id)+(c_{p-1}\tensor\id)(c_{p-1}\tensor\id)^*$ becomes a matrix over
$\Z\pi\subset\NeumannN\pi$ in this way, which acts on $l^2(\pi)^d$ via
multiplication from the left. Here
$\NeumannN\pi$ is the group von Neumann algebra with its natural trace 
$\tr_\pi$, defined as follows:
\begin{definition}
  For $\Delta=(a_{ij})\in M(d\times d,\NeumannN\pi)$ set
  \begin{equation*}
    \tr_{\pi}(\Delta):=\sum_i \tr_{\pi}(a_{ii})
  \end{equation*}
  where $\tr_{\pi}(a):=a_1=(a,1)_{l^2(\pi)}$ is the
  coefficient of the trivial group element if
  $a=\sum_{g\in\pi}\lambda_g g\in\NeumannN\pi\subset l^2(\pi)$.
\end{definition}

 Particularly important are the spectral density functions
\begin{equation}
  \label{specdens}
  F_p(\lambda):= F_{\Delta_p}(\lambda):= \tr_\pi
  \chi_{[0,\lambda]}(\Delta_p) .
\end{equation}
The $L^2$-Betti numbers
are defined as
\begin{equation*}
  b_p^{(2)}(X):= b_p^{(2)}(\Delta_p):=F_{\Delta_p}(0)=\dim_\pi(\ker(\Delta_p)) .
\end{equation*}
These are invariants of the homotopy type of $X$.

Another important invariant is the \emph{regularized determinant}:
\begin{definition}\label{defdet}
For a positive and self-adjoint operator $\Delta\in M(d\times d,\NeumannN)$ over a finite von
Neumann algebra $\NeumannN$ with spectral density function $F_\Delta$ define
\begin{equation*}
  \ln\Det_\NeumannN(\Delta) :=
  \begin{cases}
    \int_{0^+}^\infty \ln(\lambda)d F_\Delta(\lambda); &\text{if the
      integral converges}\\
    -\infty;& \text{otherwise}
  \end{cases}
\end{equation*}
Sometimes, this regularized
determinant is called \emph{Fuglede-Kadison determinant}.
\end{definition}

This gives rise to the definition:
\begin{definition}
  A self-adjoint operator $\Delta$ as above is said to be of
  \emph{$\pi$-determinant} class if and only if
\begin{equation*}
  \int_{0^+}^1{\ln\lambda}\,d F_\Delta(\lambda) > - \infty.
\end{equation*}
The space $X$ is said to be of \emph{$\pi$-determinant class} if the Laplacian
$\Delta_p$ is of determinant class for every $p$.
\end{definition}

\begin{conjecture}\label{detclass}
  Every finite $CW$-complex is of determinant class.
\end{conjecture}

If $X$ is of $\pi$-determinant class and all $L^2$-Betti numbers are
zero, then we can define its (additive) $L^2$-Reidemeister torsion
\begin{equation*}
  T^{(2)}(X):= \sum_p (-1)^p p \ln\Det_\pi\Delta_p .
\end{equation*}
Burghelea et.al.~\cite{Burghelea-Friedlander-Kappeler-McDonald(1996a)}
show that $L^2$-Reidemeister torsion is equal to $L^2$-analytical
torsion (for closed manifolds) and therefore is a generalization of
the volume of a hyperbolic manifold.

 L\"uck \cite[1.4]{Lueck(1994a)} shows that this torsion is
an invariant of the simple homotopy type of $X$. He conjectures
\begin{conjecture}
  \label{hominv}
  $L^2$-Reidemeister torsion is a homotopy invariant,
\end{conjecture}
and  proves the following theorem \cite[1.4]{Lueck(1994a)}:
\begin{theorem}
  The $L^2$-Reidemeister torsion is a homotopy invariant of
  $L^2$-acyclic finite $CW$-complexes with fundamental group $\pi$ if
  and only if
  for every (invertible!) $A\in Gl(d\times d,\Z\pi)$ the regularized
  determinant is zero:
  $\ln\Det_\pi(A^*A)=0$.
\end{theorem}
\begin{remark}
  In fact, $\ln\det_\pi(\cdot^*\cdot)$ factors through the Whitehead group of
  $\pi$. The corresponding homomorphism is denoted $\Phi$ by L\"uck,
  but we will write $\ln\det_\pi$ for the map on $Wh(\pi)$ as well.
\end{remark}

\begin{remark}
    Mathai and Rothenberg \cite[2.5]{Mathai-Rothenberg(1997)}
  extend the
  study of the $L^2$-determinants from the $L^2$-acyclic to the
  general case, dealing with determinant lines instead of complex numbers.
  Without any difficulty, this could be done in our more general
  situation as well. Because this would only complicate the notation
  and seems not to be of particular importance, we will not carry this
  out.
\end{remark}

Suppose that $\pi$ is residual, i.e.~it contains a nested sequence of normal subgroups
$\pi=\pi_1\supset\pi_2\supset\dots$ such that $\bigcap_i
\pi_i=\{1\}$. Then we construct the corresponding coverings $X_i$ of $X$
with fundamental group $\pi_i$. $L^2$-invariants are defined for
arbitrary normal coverings (the relevant von Neumann algebra is the
one of the group of deck transformations). In the given situation we
conjecture
\begin{conjecture}
  \label{resappr} In the situation just described,
  for every $p$, the $L^2$-Betti numbers $b_p^{(2)}(X_i)$ converge as
  $i\to\infty$ and
  \begin{equation*}
    \lim_{i\to\infty} b^{(2)}_p(X_i) = b^{(2)}_p(X) .
  \end{equation*}
\end{conjecture}

The projections $\pi\to\pi/\pi_i$ induce maps $p_i:M(d\times d,\Z\pi)\to
M(d\times d,\Z\pi/\pi_i)$, and the Laplacian on $X_i$ (considered as
such a matrix) is just the image of the Laplacian on $\tilde X$. Therefore,
Conjecture \ref{resappr} follows from the following conjecture about
such matrices (and, in fact is equivalent as follows from the proof of
\cite[2.2]{Lueck(1997)}).
\begin{conjecture}
  \label{matresappr} 
  For $A\in M(d\times d,\Z\pi)$ set $A_i:=p_i(A)$. Then
  \begin{equation*}
    \lim_{i\to\infty}\dim_{\pi/\pi_i}(\ker A_i)=\dim_\pi(\ker A) .
  \end{equation*}
\end{conjecture}

The Conjectures \ref{detclass}, \ref{hominv} and \ref{resappr} are
proven for residually finite  groups by L\"uck
\cite[0.5]{Lueck(1994c)}. Using the
ideas of L\"uck in a different context, Dodziuk and Mathai prove
Conjecture \ref{detclass} for amenable $\pi$
\cite[0.2]{Dodziuk-Mathai(1996a)}. They also
establish an approximation theorem for $L^2$-Betti numbers
of a slightly different type in this case.

Michael Farber \cite{Farber(1997)} generalises L\"uck's results to
sequences of finite dimensional coefficients, which converge to
$l^2(\pi)$ (one can interpret the covering cohomology as cohomology
with coefficients in $l^2(\pi/\pi_i)$). In this sense, we are
interested in special but infinite dimensional coefficient systems.

The aim of this paper is to extend the results of L\"uck and
Dodziuk/Mathai to the following larger class of groups (Conjectures
\ref{resappr} and \ref{matresappr} have to be extended and  modified suitably).

\begin{definition}
  Let $\RAgroups$ be the smallest class of groups which contains the trivial
  group and is closed under the following processes:
  \begin{itemize}
   \item If $U\subgroup \pi$ is any subgroup such that $U\in\RAgroups$
    and the discrete homogeneous space $\pi/U$ admits a
    $\pi$-invariant metric which makes it to an amenable discrete
    metric space, then $\pi\in\RAgroups$. (For our notion of
    amenability compare Section \ref{sec:am}. The most important
    example is:  $U$  normal and $\pi/U$ is an amenable group.)
  \item If $\pi=\dirlim_{i\in I}\pi_i$ is the direct limit of a
    directed system of groups $\pi_i\in \RAgroups$,
    then  $\pi\in\RAgroups$, too.
  \item If $\pi=\invlim_{i\in I}\pi_i$ is the inverse limit of a
    directed system of groups $\pi_i\in \RAgroups$, then 
    $\pi\in\RAgroups$, too.
  \item The class $\RAgroups$ is closed under taking subgroups.
  \end{itemize}
\end{definition}

\begin{remark}
  It follows immediately from the definition that $\RAgroups$ contains
  all amenable groups, is
  closed under directed unions and is residually closed. More details
  can be found in section \ref{sec:RAprop}.
\end{remark}

The main theorem of the paper is the following:
\begin{theorem}\label{main}
  Suppose $\pi$ belongs to the class $\RAgroups$. Then for every
  CW-complex with fundamental group $\pi$ which has finitely many
  cells in each dimension, the Conjectures \ref{detclass} and
  \ref{hominv} are true.  Approximation results
  which generalize   \ref{resappr} and \ref{matresappr} are valid
  under the condition that all the occuring groups belong
  to $\RAgroups$.
\end{theorem}

\begin{remark}
  Clair
  proves in \cite{Clair(1997)}
  the Conjectures \ref{detclass}, \ref{hominv} and \ref{resappr} for
  the large class residually amenable 
  fundamental groups, using \cite[2.5]{Mathai-Rothenberg(1997)}.

  Our general principle is to use the methods of L\"uck (and
  Dodziuk, Mathai, Ro\-then\-berg) and to check carefully which is the most
  general situation they apply to.
\end{remark}

\begin{remark}
  So far, no example of a countable group which does not belong to the
  class $\RAgroups$ has been constructed. Good candidates for such examples
  are finitely generated simple groups which are not amenable,
  e.g.\ groups containing a free group with two generators.

  On the other hand, we can give no example of a non-residually
  amenable group which belongs to $\RAgroups$, either. In any case, our
  description of $\RAgroups$ leads easily to many properties like
  closedness under direct sums and free products, which (if true at
  all) are probably much harder to establish for the class of
  residually amenable groups.
\end{remark}

In fact, we prove a little bit more than Theorem \ref{main}. Namely,
we show that the relevant
properties  are stable under the
 operations characterizing the class $\RAgroups$. We use the
 following definitions:
\begin{definition}
  Let $C$ be any property of discrete groups. It is said to be
  \begin{itemize}
  \item \emph{stable under direct/inverse limits} if $C$ is true for
    $\pi$ whenever $\pi$ is a direct/inverse limit of a directed
    system of groups which have property $C$.
  \item \emph{subgroup stable} if every subgroup $U\subgroup \pi$ of a
    group with property $C$ shares this property, too.
                                % stable under finite extensions
  \item \emph{stable under amenable extensions} if $\pi$ has property
    $C$ whenever it contains
    a subgroup $U$ with property $C$ such that the homogeneous space
    $\pi/U$ is amenable.
  \end{itemize}
\end{definition}

The properties we have in mind are listed in the following definition:
\begin{definition}\label{prop_list}
  Let $\pi$ be a discrete group. We say
  \begin{itemize}
  \item $\pi$ is of determinant class, if $\Delta$ is of
    $\pi$-determinant class $\forall \Delta\in M(d\times
    d,\integers\pi)$ which is positive and self-adjoint;
  \item $\pi$ has \emph{semi-integral determinant} if $\ln\det_\pi(\Delta)\ge
    0$ for every $\Delta\in M(d\times d,\Z\pi)$ which is positive
    self-adjoint. In particular every such $\Delta$ is of
    $\pi$-determinant 
    class, i.e.~$\pi$ itself is of determinant class;
  \item $\pi$ has \emph{Whitehead-trivial determinant} if
    $\ln\det_\pi(A^*A)=0$ $\forall A\in Wh(\pi)$.
  % spectral sublog.
  \end{itemize}
\end{definition}

\begin{theorem}\label{hominvstab}
  Whitehead-trivial determinant is stable under
  direct and inverse limits and is subgroup stable.
\end{theorem}

\begin{remark}
  In the light of this theorem, we can use the fact that the
  Fuglede-Kadison determinant must
  be trivial on trivial
  Whitehead groups, e.g.\ for every torsion free discrete and
  cocompact subgroup of a Lie group with finitely many components
  \cite[2.1]{Farrell-Jones(1993a)}. Waldhausen shows that the
  Whitehead group is trivial for another class of groups,
  including torsion free one-relator groups and many fundamental
  groups of 3-manifolds
  \cite[17.5]{Waldhausen(1978a)}.

 The validity of the isomorphism conjecture of Farell
  and Jones would imply that
  the Whitehead group is trivial if $\pi$ is torsion free.

  We can define the class $\RAgroups'$ as the smallest
  class of groups which contains $\RAgroups$, in addition any class of
  groups whose Whitehead group is (known to be) trivial, and which is
  closed under
  taking subgroups, direct and inverse limits. Then every group in
  $\RAgroups'$ has Whitehead trivial determinant, i.e.~$L^2$-torsion
  is a homotopy invariant for $L^2$-acyclic finite $CW$-complexes with
  such a fundamental group.
\end{remark}

\begin{theorem}\label{semiint_stab}
  The property ``semi-integral determinant'' is stable under direct and inverse limits,  subgroup stable  and  stable under  amenable extensions.
%vorher noch:  and ``spectral sublogarithmicity''  are
\end{theorem}

\begin{theorem}\label{T20}
  If for a group $\pi$ and $\forall A\in Wh(\pi)$ we have
  $\ln\Det_\pi(A^*A)\ge 0$, then
  the Fuglede-Kadison determinant is trivial on $Wh(\pi)$.

  In particular, semi-integral determinant implies Whitehead-trivial
  determinant.
\end{theorem}
\begin{proof}
  $A\in Wh(\pi)$ implies $A$ has an inverse $B\in Wh(\pi)$. Now by
  \cite[4.2]{Lueck(1994a)}
  \begin{equation*}
    0=\ln\Det_\pi(\id)=\ln\Det_\pi((AB)^*AB)
    =\underbrace{\ln\Det_\pi(A^*A)}_{\ge
      0}+\underbrace{\ln\Det_\pi(B^*B)}_{\ge 0}
  \end{equation*}
  and the statement follows.
\end{proof}

It follows from the induction principle \ref{ind_princ} and the fact
that the trivial group has semi-integral determinant (Lemma
\ref{trivgroup}) that
\ref{semiint_stab} and \ref{T20} imply the first part of our main
theorem \ref{main}.

\section{Properties of the class $\RAgroups$}\label{sec:RAprop}

We proceed with a more precise definition of the class $\RAgroups$,
similar to the description Linnel gives of his
class $\mathcal{C}$ in \cite[p.~570]{Linnel(1993)}.

\begin{definition}
  For each ordinal $\alpha$, we define a class of groups
  $\RAgroups_\alpha$ inductively.
  \begin{itemize}
  \item $\RAgroups_0$ consists  of the trivial groups.
  \item If $\alpha$ is a limit ordinal, then $\RAgroups_\alpha$ is the
    union of $\RAgroups_\beta$ with $\beta<\alpha$
  \item If $\alpha$ has a predecessor $\alpha-1$, then
    $\RAgroups_\alpha$ consists of all groups which
    \begin{itemize}
    \item are subgroups of groups of $\RAgroups_{\alpha-1}$
    \item contain a subgroup $U$ which belongs to
      $\RAgroups_{\alpha-1}$ such that the quotient space is an amenable
      homogeneous space
    \item are direct or inverse limits of directed systems of groups
      in $\RAgroups_{\alpha-1}$.
    \end{itemize}
  \end{itemize}
  By definition, a group is in $\RAgroups$ if it belongs to $\RAgroups_\alpha$
  for some ordinal $\alpha$.
\end{definition}

\forget{
\begin{remark}
  Every class of groups which is closed under taking subgroups, direct
  and inverse limits (of directed systems) and amenable extensions in
  the above sense must contain $\RAgroups$.

  On the other hand, if
  $\pi\in\RAgroups$, say $\pi\in\RAgroups_\alpha$, and $U\subgroup \pi$,
  then $U\in\RAgroups_{\alpha+1}$, similarly for amenable
  extensions. If $\pi$ is the limit of $(\pi_i)_{i\in I}$ with
  $\pi_i\in\RAgroups_{\alpha_i}$, then let $\alpha$ be the smallest
  ordinal larger than each of the $\alpha_i$ and it follows
  $\pi\in\RAgroups_\alpha$. 
\end{remark}
}

The class $\RAgroups$ is defined by (transfinite)
  induction. Therefore, properties of the groups in $\RAgroups$ can be
  proven by induction, too. More precisely, the following induction
  principle is valid:
  \begin{proposition}\label{ind_princ}
    Suppose a property $C$ of groups is shared by the trivial group, and
    the following is true:
    \begin{itemize}
    \item  whenever $K$ has property $C$ and $K\subgroup \pi$
      with $\pi/K$ an amenable
      homogenous space, then  $\pi$ has property $C$ as well;
    \item whenever $\pi$ is a direct or inverse limit of a directed
      system of groups with the property $C$, then $\pi$ has
      property $C$
    \item property $C$ is inherited by subgroups.
    \end{itemize}
    Then property $C$ is shared by all groups in
    the class $\RAgroups$.
  \end{proposition}
  \begin{proof}
    The proof of the induction principle is done by transfinite
    induction.

    By assumption, $C$ holds for $\RAgroups_0$. We have to establish
    $C$ for every group in $\RAgroups_\alpha$, granted its validity
    for groups in $\RAgroups_\beta$ for all $\beta<\alpha$. If
    $\alpha$ is a limit ordinal, this is trivial. If $\alpha$ has a
    predecessor $\alpha-1$, the assumptions just match the definition
    of $\RAgroups_\alpha$, so the statement follows.
  \end{proof}

Now, we study the properties of the class $\RAgroups$.
\begin{proposition}
  The class $\RAgroups$ is closed under directed unions.
\end{proposition}
\begin{proof}
    A directed union is a special case of a directed direct limit.
\end{proof}

\begin{proposition}
  The class $\RAgroups$ is
  residually closed. This means that if $\pi$
  contains a nested sequence of normal subgroups
  $\pi_1\supset\pi_2\supset\dots$ with trivial intersection and if
  $\pi/\pi_i\in\RAgroups$ $\forall i$ then also $\pi\in\RAgroups$.
\end{proposition}
\begin{proof}
  The inverse system of groups $\pi/\pi_i$ has some inverse
  limit $G$. The system of maps $\pi\to\pi/\pi_i$ induces a
  homomorphism $\pi\to G$. If $g\in\pi$ is mapped to $1\in G$, then
  $g$ has to be mapped to $1\in\pi/\pi_i$ $\forall i$,
  i.e.~$g\in\bigcap_i \pi_i=\{1\}$. As a directed limit, $G\in
  \RAgroups$, and as a subgroup of $G$, also $\pi\in\RAgroups$, as well.
\end{proof}

\begin{theorem}
  If $U$ belongs to $\RAgroups$ and $i:U \to U$ is any group
  homomorphism, then the ``mapping torus''-extension of $U$ with respect to $i$
  \begin{equation*}
    \pi = \generate{u\in U,t| t^{-1}ut = i(u), u\cdot v= (uv); \forall u,v\in U}
  \end{equation*}
  also belongs to $\RAgroups$ (if $i$ is injective, this is a special
  example of an HNN-extension).
\end{theorem}
\begin{proof}
  There is a canonical projection $\pi\to \Z$ sending $u\in U$ to $0$
  and $t$ to $1$. Denote its kernel by $K$. We will show that $K$
  belongs to $\RAgroups$, then so does  $\pi$  because it is an
  extension of $K$ with amenable quotient $\Z$.

  Now $K$ is the direct limit of the sequence 
  \begin{equation*}
    U\xrightarrow{i} U \xrightarrow{i} U \xrightarrow {i} U \dots
  \end{equation*}
  and belongs to $\RAgroups$, which is closed under taking direct limits.
\end{proof}

\begin{remark}
  Although such a mapping-torus extension of a finitely presented
  group is finitely
  presented again, the kernel $K$ we used in the proof may very well
  not even admit a finite set of generators. This is one instance
  where it is useful to allow arbitrary groups, even if one is only
  interested in fundamental groups of finite $CW$-complexes.
\end{remark}

For the next property, we use the induction principle. 
\begin{proposition}\label{freeprod}
  $\RAgroups$ is closed under forming
  \begin{enumerate}
  \item  direct sums and direct products;
  \item free products.
\end{enumerate}
\end{proposition}
\begin{proof}
We have to check the conditions for the induction principle. Fix  $\pi\in\RAgroups$.
\begin{enumerate}
\item If $U\subgroup G$ then $U\times \pi\subgroup G\times\pi$. If
  $G\times \pi\in\RAgroups$, the same is true for
  $U\times\pi$. If $G$ is the (direct or inverse) limit of the
  directed system of groups $G_i$, then $G\times \pi$ is the limit of
  the system $G_i\times \pi$ (compare Lemma \ref{lim1} or
  \ref{lim2}). If (by assumption) $G_i\times
  \pi\in\RAgroups$, then $G\times\pi\in\RAgroups$.
  Finally, if $U\subgroup G$, and $G/U$ is amenable, then $U\times
  \pi\subgroup G\times \pi$ with the same amenable quotient.
  Therefore, $U\times\pi \in\RAgroups$ $\implies$ $G\times\pi\in\RAgroups$.
\item 
  For the free products, the first step is to prove: $*_{i\in
    I}\Z/4\in\RAgroups$ for every index set $I$. This is the direct limit of
  finite free
  products of copies of $\Z/4$, therefore we have to prove the
  statement for finite $I$. Now $*_{i=1}^n\Z/4$ is a subgroup of
  $\Z/4*\Z/4$ (contained in the kernel of the projection onto one
  factor), and $\Z/4*\Z/4$ is virtually free, i.e.~an extension of a
  residually finite (the free group) with an amenable (finite) group.
  Therefore  $\integers/4*\integers/4$ belongs to $\RAgroups$.
  
  Next we show that $\pi*(*_{j\in J}\integers/4)\in\RAgroups$ for every
  set $J$.We prove this using the induction principle. For $\pi=1$ this is the
  conclusion of the first step. If $\pi$ is a limit of $(\pi_i)_{\in
  I}$, or a subgroup of $G$, then $\pi* (*_{j\in J}\Z/4)$
  is a subgroup of the limit of $\pi_i*(*_{j\in J}\Z/4)$ (compare
  Lemma \ref{lim1} and Lemma \ref{lim2}) or a subgroup of
  $G*(*_{j\in J}\Z/4)$, and we can apply that 
  $\RAgroups$ is subgroup closed. If $U\subgroup \pi$ and $\pi/U$ is
  amenable,  $\pi*(*_{j\in J}\integers/4)$ acts on $\pi/U$. We
  get a new point stabilizer, which is isomorphic to the free product
  of $U$ with $*_{G/U}(*_{i\in I}\integers/4)$. Fortunately, the induction
  hypothesis applies with the free product of an arbitrary number of
  copies of $\integers/4$.
  
  As the next step we show that  $*_{i\in I}\pi\in\RAgroups$. This follows
  (as $*_{i\in I}\pi$ is a direct limit) from the
  corresponding statement for $I$ finite, and these are
  subgroups of $\pi*\Z/4$, contained in the kernel of the projection
  onto $\Z/4$.
  
  $\pi_1*\pi_2$ is contained in $(\pi_1\times \pi_2)*(\pi_1\times
  \pi_2)$, and the general statement follows by induction and taking
  limits. \qquad\qed
\end{enumerate}
\renewcommand{\qed}{}
\end{proof}

In the proof of  Proposition \ref{freeprod} we have used the following
two lemmas.
\begin{lemma}\label{lim1}
  If $\pi$ is the direct limit of a system of groups $\pi_i$ and $G$
  is any group, then $\pi*G$ is the direct limit of $\pi_i*G$, and
  $\pi\times G$ is the direct limit of $\pi_i\times G$.
\end{lemma}
\begin{proof}
  There are obvious maps from $\pi_i*G$ to $\pi*G$ and from
  $\pi_i\times G$ to $\pi\times G$. 
  
  Suppose one has consistent maps from $\pi_i*G$ (or $\pi_i\times G$)
  to some group
  $X$. Since $\pi_i$ and $G$ both are subgroups of $\pi_i*G$ (or of
  $\pi_i\times G)$ , this
  means that we have a consistent family of maps on $\pi_i$
  multiplied with a fixed map on $G$. Therefore (from the properties
  of products) there exists exactly one map from $\pi*G$ (or
  $\pi\times G$) to $X$ 
  making all the diagrams commutative (for the commutative product
  note that the union of the images of the $\pi_i$ in $X$
  commutes with the image of $G$).
\end{proof}

\begin{lemma}\label{lim2}
  If $\pi$ is an inverse limit of a system of groups $\pi_i$ and $G$
  is any group, then $\pi*G$ is contained in the inverse limit $X$ of
  $\pi_i*G$.

  The inverse limite of $\pi_i\times G$ is $\pi\times G$.
\end{lemma}
\begin{proof}
  First, we look at the free products:
  
  We have a consistent family of homomorphisms from $\pi*G$ to $\pi_i*G$,
  therefore a homomorphism from $\pi*G$ to $X$. An element
  $x=p_1g_1\dots p_n g_n\in \pi*G$ is in the kernel of this homomorphism
  iff it is mapped to $1\in\pi_i*G$ for every $i\in I$. This can not
  happen if $1\ne x\in\pi$. It remains to check the case $1\ne g_1\in
  G$. We may assume that $g_2\ne 1$ iff $p_2\ne 1$. If
  $\phi_i:\pi\to\pi_i$ is the natural homomorphism, then $x$ is mapped
  to $\phi_i(p_1)g_1\phi_i(p_2)g_2\dots g_n\in\pi_i*G$. If this is
  trivial, but $g_1\ne 1$, necessarily $\phi_i(p_2)=1$ $\forall
  i$. This implies $p_2=1$, i.e.~$x=p_1g_1$, since we wrote $x$ in
  normal form. But then $\phi_i*\id(x)\ne 1$ $\forall i$, and the
  kernel of the map to $X$ is trivial, as required.

  For the commutative product, let $X$ be a group together with a
  consistent family of morphisms to $\pi_i\times G$. These have the
  form $x\mapsto (\phi_i(x),f_i(x))$. Composition with the projections
  to $\pi_i$ or to $G$ shows that $\pi_i$ is a consistent family of
  morphisms to $\pi_i$, and $f=f_i:X\to G$ all coincide. Let $\phi:X\to
  \pi$ be the limit. Then $\phi\times f:X\to \pi\times G$ is a unique
  homomorphism which makes all relevant diagrams commutative.
\end{proof}

\section{Passage to subgroups}

Suppose $U\subset\pi$ is a subgroup of a discrete group. A positive
self-adjoint matrix $A\in M(d\times d,\Z U)$ can also be considered as
a matrix over $\Z\pi$. Denote the operators with $A_U$ and $A_\pi$,
respecitively. Recall the following well known fact:
\begin{proposition}\label{inddens}
  The spectral density functions of $A_U$  and $A_\pi$ coincide.
\end{proposition}
\begin{proof}
  Choose a set of representatives $\{g_i\}_{i\in I}$ with $0\in I$ and
  $g_0=1$, to write
  $\pi=\disjointunion_{i\in I} Ug_i$. Then 
  \begin{equation*}
    l^2(\pi)^d =\bigoplus_{i\in I} l^2(U)^d g_i .
  \end{equation*}
With respect to this splitting, the action of $A_\pi$ on $l^2(\pi)$ is
diagonal and, restricted to each of the summands $l^2(U)^d g_i$, is
multiplication by $A_u$ from the left. It follows that every spectral
projection $\chi_{[0,\lambda]}(A_\pi)$ is diagonal with
$\chi_{[0,\lambda]}(A_U)$ on the diagonal. Then
\begin{equation*}
  \begin{split}
    F_{A_\pi}(\lambda) =& \sum_{k=1}^d
    \innerprod{\chi_{[0,\lambda]}(A_\pi)e_k^\pi,e_k^\pi}
    = \sum_{k=}^d \innerprod{\chi_{[0,\lambda]}(A_\pi) e_k^U\cdot
      1,e_k^U\cdot q}\\
    =& \sum_{k=1}^d \innerprod{\chi_{[0,\lambda]}(A_U)e_k^U,e_k^U} =
    F_{A_U}(\lambda) . \qquad\qed
  \end{split}
\end{equation*}
\renewcommand{\qed}{}
\end{proof}

\begin{corollary}
  The properties of Definition \ref{prop_list} are inherited by subgroups.
\end{corollary}
In particular, we have proven the subgroup part of Theorems
\ref{hominvstab} and \ref{semiint_stab}.

\section{Amenable extensions}\label{sec:am}

\begin{definition}
  A discrete homogeneous space $\pi/U$ is called \emph{amenable}, if on $\pi/U$
  we find a $\pi$-invariant integer-valued metric $d:\pi/U\times \pi/U\to\N$
  such that 
  \begin{itemize}
  \item sets of finite diameter are finite 
  \item for every $K>0$,
    $\epsilon>0$ there is a finite subset $X\subset \pi/U$ with
  \begin{equation*}
    \abs{ N_K(X)} \le \epsilon\abs{ X}
  \end{equation*}
  where $N_K(X):=\{x\in \pi/U;\; d(x,X)\le K\text{ and }d(x,\pi/U-X)\le
  K\}$ is the $K$-neighborhood
  of the boundary of  $X$.
\end{itemize}
 A nested sequence of finite subsets $K_1\subset K_2\subset \dots$ is called
  an \emph{amenable exhaustion} of $\pi/U$ if $\bigcup K_n=\pi/U$ and if
  $\forall K>0$ and $\epsilon>0$  we find $N\in\N$ so
  that $\abs{ N_K(K_i)}\le \epsilon\abs{ K_i}$ $\forall i\ge N$.
\end{definition}

\begin{lemma}
  Every amenable homogeneous space $\pi/U$ admits an amenable exhaustion.
\end{lemma}
\begin{proof}
  By assumption for $n,K\in\N$ we find $X_{n,K}$ with $\abs{
  N_K(X_{n,K})}\le \frac{1}{n}\abs{ X_{n,K}}$. Fix some base point in
  $\pi/U$. Since $\pi$ acts transitively on $\pi/U$ and the metric is
  $\pi$-invariant, we may assume after translation that the base point
  is contained in each of the $X_{n,K}$. Now we construct 
  the exhaustion $E_i$. Set $E_1:= X_{1,1}$ inductively. For the
  induction suppose
  $E_1,\dots, E_n$ are constructed with $\abs{ N_{k}(E_k)}\le \frac{1}{k}\abs{
  E_i}$ for $k=1,\dots,n$. Suppose $E_n$ has diameter $\le d\in\N$ with
  $2d\ge n+1$. Set $E_{n+1}:= E_n\cup X_{n+1,2d}$. Then
  $N_{n+1}(E_{n+1})\subset N_{2d}(E_{n+1})$. On the
  other hand, by the triangle inequality $N_d(E_n)\subset
  X_{n+1,2d}$ and therefore
  $N_{2d}(E_{n+1})=N_{2d}(X_{n+1,2d})$, and it follows
  \begin{equation*}
    \abs{ N_{n+1}(E_{n+1})}\le \abs{ N_{2d}(X_{n+1,2d})}\le \frac{1}{n+1}
    \abs{ X_{n+1,2d}}\le \frac{1}{n+1}\abs{ E_{n+1}} .
  \end{equation*}
  The claim follows.
\end{proof}

\begin{example}
  If $U$ is a normal subgroup and $\pi/U$ is an amenable group, it is an
  amenable homogeneous space, too.
\end{example} 

\begin{definition}\label{amsit}
  Suppose $\pi$ is a group with subgroup $U$ and amenable quotient
  $\pi/U$. Choose an amenable cover $X_1\subset X_2\subset \dots\subset
  \pi/U$. For 
  $B\in M(d\times d,\NeumannN U)$ set
  \begin{equation*}
    \tr_m(B):=\frac{1}{\abs{X_m}}\tr_U(B).
  \end{equation*}
  For
 $\Delta\in M(d\times d,\NeumannN \pi)$ positive and self-adjoint  set
 $\Delta_n:= P_m\Delta P_m$
  where $P_m=\diag(p_m)$ with $p_m\in \boundedops(l^2(\pi))$ is given by
  projection onto the closed subspace generated by the inverse image
  of $X_m$. Then $\Delta_m$ no longer belongs to $\NeumannN \pi$ but
  still to $\NeumannN U$ and we define (by slight abuse of notation)
  \begin{equation*}
    \begin{split}
      F_{\Delta_m}(\lambda) &:=
      \tr_m(\chi_{[0,\lambda]}(\Delta_m)) ,\\
      \ln\Det_U(\Delta_m) &:=\int_{0^+}^\infty
      \ln(\lambda)\;dF_{\Delta_m}(\lambda)\qquad\text{using the new
        $F_{\Delta_m}$}.
     \end{split} 
  \end{equation*}
  Here $\Delta_m$ is considered as operator on the image of $P_m$.
  This subspace is $\NeumannN U$-isomorphic to $l^2(U)^{d\abs{X_m}}$.

  Note that there are two meanings of $F_{\Delta_m}(\lambda)$ and
  $\ln\det_U(\Delta_m)$  (using
  either $\tr_U$ or $\tr_m$), but in
  the amenable case we will always use the variant where we divide by
  the volume of the sets $X_m$.
\end{definition}

The following is one of the key lemmas which make our (respectively
L\"uck's) method work:
\begin{lemma}\label{ambound}
  In the situation above, we find $K\in\reals$ independent of $m$, so
  that 
  \begin{equation*}
    \norm{\Delta}\le K\qquad\text{and}\qquad\norm{\Delta_m}\le K\qquad\forall m\in\naturals.
  \end{equation*}
\end{lemma}
\begin{proof}
  This is an immediate consequence of the fact that $\norm{P}\le1$ for every
  projection $P$  and $\Delta_m=P_m\Delta P_m$ with
  projections $P_m$.
\end{proof}

We now establish the second key lemma. It
generalizes a
corresponding result of Dodziuk/Mathai
\cite[2.3]{Dodziuk-Mathai(1996a)} where $U$ is trivial. We need the
result only for matrices over $\C\pi$, but for possible other
applications we proof a more general statement here.

\begin{lemma}\label{amappr}
     Let $p(x)\in\C[x]$ be a polynomial. Suppose $\Delta\in M(d\times
     d,\NeumannN\pi)$. Then 
  \begin{equation*}
  \tr_\pi
    p(\Delta)=\lim_{m\to\infty}\tr_m p(\Delta_m).
  \end{equation*} 
\end{lemma}
  \begin{proof}
By linearity it suffices to prove the statement for the monomials
$x^N$, $N\in\naturals$.

Pull the
    metric on $\pi/U$ back
    to $\pi$ to get some semimetric on $\pi$. Denote the inverse image of
    $X_k$ in $\pi$ with $X_k'$.

We have to compare $(\Delta^N g e_k,g e_k)$ and
$(\Delta_m^N g e_k,g e_k)$ for $g=g_i\in X_m'$, in particular for
those $g_i$ with $B_a(g_i)\subset X_m'$. Of course, we don't find $a\in\reals$
such that the difference is zero. However, we will show that (for
fixed $N$) we can find $a$ such that the difference is sufficiently small.

First observe that $P_m g e_k=g e_k$ if $g\in X_m'$, and (since $P_m$ is self-adjoint)
\begin{equation*}
  ( (P_m\Delta P_m)^N g e_k,g e_k) =  (\Delta P_m\Delta\dots P_m\Delta
  g e_k,g e_k) .
\end{equation*}
Now the following sum is a telescope and therefore
  \begin{multline}
    \Delta P_m\Delta \dots P_m\Delta =\\
     \Delta^N -
    \Delta(1-P_m)\Delta^{N-1} -\Delta
    P_m\Delta(1-P_m)\Delta^{N-2}-\dots - \Delta P_m\dots \Delta
    (1-P_m)\Delta .
\end{multline}
It follows for $g\in X_m'$
\begin{equation*}
  \begin{split}
    \abs{ (\Delta^Ng e_k,g e_k)-(\Delta_m^Ng e_k,g e_k)}
 &\le \sum_{i=1}^{N-1} \abs{(
      (1-P_m)\Delta^i g e_k, (\Delta^*P_m)^{N-i}g e_k)}\\
    & \le \sum_{i=1}^{N-1}\abs{(1-P_m)\Delta^i g
      e_k}\cdot\norm{\Delta^*}^{N-i} .
  \end{split}
\end{equation*}
Here we used the fact that the norm of a nontrivial projector is $1$
and $\abs{g e_k}=1$.

%\Kommentar{Schrott, der ueberig war:}
% &\le
%        \frac{1}{\abs{X_m}}\sum_{k=1}^d \sum_{i\in
%          N_{ra}(X_m)}\abs{(p(\Delta)g_i e_k,g_i e_k)}
%        +\abs{(p(\Delta_m)g_i e_k,g_i e_k)}\\
%        & \le \frac{\abs{ N_{ar}(X_m)}}-(\Delta_m^Ng e_k,g e_k)\\
Fix $\epsilon>0$. For $i=1,\dots,N-1$ and $ k=1,\dots,d$ we have $\Delta^i g e_k\in l^2(\pi)^d$. It
follows that we find $R>0$ so that
\begin{equation}\label{eps_ie}
  \abs{ (1-P_{B_R(g)})\Delta^i g e_k}\le\epsilon
\end{equation}
where $P_{B_R(g)}$ is the projector onto the closed subspace
spanned by the elements in $\bigcup_{k=1}^d B_R(g)e_k$. Since $\Delta$
and the semimetric are $\pi$-invariant, this holds for every $g\in
\pi$ with $R$
independent of $g$. If the range of $P_{B_R(g)}$ is contained in the
range of $P_m$, i.e.\ if $B_R(g)\in X_m'$ then \eqref{eps_ie} implies
\begin{equation*}
  \abs{ (1-P_m)\Delta^i g e_k}\le\epsilon
\end{equation*}
(since we have even more trivial Fourier coefficients in the standard
orthonormal base coming from $\pi$ of
$l^2(\pi)^d$). Taken together, we get
\begin{equation*}
  \begin{split}
    &\abs{\tr_\pi \Delta^N-\tr_m \Delta_m^N}  \le 
     \frac{1}{\abs{X_m}}\sum_{k=1}^d\sum_{i\in X_m}\abs{(\Delta^N g_i e_k,g_i
      e_k)-(
      (\Delta_m)^N g_i e_k,g_i e_k)} \\
     & \le \frac{1}{\abs{X_m}} \sum_{k=1}^d\sum_{i\in X_m}\sum_{j=1}^{N-1}
    \abs{(1-P_m)\Delta^j g_i e_k}\cdot\norm{\Delta^*}^{N-j}\\
     & \le  \frac{1}{\abs{X_m}}\sum_{k=1}^d\sum_{i\in X_m-N_R(X_m)}\sum_{j=1}^{N-1}\underbrace{
      \abs{(1-P_m)\Delta^j g_i
        e_k}}_{<\epsilon\text{ since $B_R(g_i)\subset
        X_m'$ by definition of $N_R(X_m)$}}\norm{\Delta^*}^{N-j}\\
    & \quad+
    \frac{1}{\abs{X_m}}\sum_{i\in N_R(X_m)}\sum_{k=1}^d\sum_{j=1}^{N-1}
    \abs{(1-P_m)\Delta^j g_i e_k}\cdot\norm{\Delta^*}^{N-j}\\
    & \le  \epsilon \underbrace{d N
      \max_{j=1,\dots,N-1}\{\norm{\Delta^*}^j\}}_{=:C_N} +
    \frac{\abs{N_R(X_m)}}{\abs{X_m}} \underbrace{d N
    \underbrace{\norm{1-P_m}}_{\le 2}\max_{j=1\dots N}\{\norm{\Delta}^j\cdot\norm{\Delta^*}^{N-j}\}}_{=:C_N'}
  \end{split}
\end{equation*}
Note that $C_N$ and $C_N'$ are independent of $m$ and
$\epsilon$. Since $X_m$ is an amenable extension, for every $R$ we
find $m_R$ so that $\frac{\abs{N_R(X_m)}}{\abs{X_m}}$ is smaller than
$\epsilon$ for every $m\ge m_R$. Since $\epsilon$ was arbitrary, the
assertion of the lemma follows.
  \end{proof}

\forget{
\begin{remark}\label{geochange}
  In our geometric setting, $\Delta$ has the shape $\Delta=A^*A+BB^*$
  with $A\in M(d\times d_1,\complexs \pi)$, $B\in M(d_2\times
  d,\complexs \pi)$. It is possible to work with $\Delta'_m=P_mA^*P_m A
  P_m + P_mB P_m B^* P_m$. The above proofs also apply to this
  situation and establish the key lemmas for this family of
  operators. These operators are closer to geometry, because they
  arise as Laplacians of subcomplexes which approximate the large
  complex whose combinatorial Laplacian is $\Delta$.
\end{remark}
}

\section{Direct and inverse limits}

\begin{remark}\label{oBdA}
In this section, we study the properties of Definition \ref{prop_list} for
direct and inverse limits.
  However, we will only deal with the apparently weaker statements that each
  of the conditions holds for every $\Delta$ of the form
  $\Delta=A^*A$. The general case is a consequence of this since we 
  can easily compare the self-adjoint $\Delta$ with
  $\Delta^2=\Delta^*\Delta$, because $F_{\Delta}(\lambda)=F_{\Delta^2}(\lambda^2)$.
\end{remark}

We describe now the situation we are dealing with in this section:
\begin{definition}\label{ressit}
  Suppose the group $\pi$ is the direct or inverse limit of a
  directed system of groups $\pi_i$, $i\in I$. The latter means that we have
  a partial ordering $<$ on $I$, and $\forall i,j\in I$ we find $k\in
  I$ with $i<k$ and $j<k$. In the case of a direct limit, let
  $p_i:\pi_i\to\pi$ the natural maps, in the
  case of an inverse limit, $p_i:\pi\to \pi_i$.

Suppose $A\in M(d\times d,\C\pi)$ is given.

  If $\pi$ is an inverse limit, let $A_i=p_i(A)$ be the image of $A$
  under the projection
  $M(d,\pi)\to M(d,\pi_i)$. Set $\Delta:=A^*A$. Then
  $\Delta_i=(A_i)^*A_i$ (this follows from the algebraic description of
  the adjoint \cite[p.~465]{Lueck(1994c)}). In particular, all of the
  operators $\Delta_i$ are positive. Define
  \begin{equation*}
    \begin{split}
      \tr_i(\Delta_i) &:=\tr_{\pi_i}(\Delta_i) .
    \end{split}
\end{equation*}
    $F_{\Delta_i}(\lambda)$ is defined using the trace on the von
    Neumann algebra $\pi_i$ ($i\in I$).

If we want to give a similar definition in the case where $\pi$ is a
direct limit, we have to make additional choices. Namely, let
$A=(a_{kl})$ with $a_{kl}=\sum_{g\in\pi}\lambda^g_{kl}g$. Then, only
finitely many of the $\lambda^g_{kl}$ are nonzero. Let $V$ be the
corresponding finite collection of $g\in\pi$. Since $\pi$ is the
direct limit of $\pi_i$ we find $j_0\in I$ such that $V\subset
p_{j_0}(\pi_{j_0})$. Choose an inverse image for each $g$ in
$\pi_{j_0}$. This gives a matrix $A_{j_0}\in M(d\times
d,\pi_{j_0})$, which is mapped to $A_i\in M(d\times d,\pi_i)$ for
$i>j_0$. Now we apply the above constructions to this net
$(A_i)_{i>j_0}$. Note that this definitely depends on the choices.

For notational convenience, we choose some $j_0\in I$ also when we deal
with an inverse limit.
\end{definition}

Now, we will establish in this situation the two key lemmas corresponding
to Lemma \ref{ambound} and Lemma \ref{amappr}.

For the first lemma, instead of working with the norm of operators, we will use
another invariant which gives an upper bound for the norm but is much
easier to read off:
\begin{definition}
  Let $\pi$ be a discrete group, $\Delta\in M(d\times
  d,\integers\pi)$. Set
  \begin{equation*}
    K(\Delta):= d^2\max_{i,j}\{\abs{a_{i,j}}_1\}\quad\text{where $\abs{\cdot}_1$
     is the $L^1$-norm on $\C\pi\subset l^1(\pi)$.}
  \end{equation*}
\end{definition}

\begin{lemma}
  \label{limbound}
  Adopt the situation of Definition \ref{ressit}. One can find
  $K\in\reals$, independent of $i$, such that
  \begin{equation*}
    \norm{A_i} \le K\quad\forall i>j_0\qquad\text{and}\qquad\norm{A}\le K .
  \end{equation*}
\end{lemma}
\begin{proof}
  L\"uck \cite[2.5]{Lueck(1994c)} shows that $\norm{A_i}\le K(A_i)$. It follows from the
  construction of $A_i$ that $K(A_i)\le K(A)$ in the case of an
  inverse limit, and $K(A_i)\le K(A_{j_0})$ in the case of a direct
  limit, with $j_0$ as above. In both cases, we obtain a uniform bound
  for $\norm{A_i}$.
\end{proof}

\begin{lemma}\label{limappr} Adopt the situation of Definition \ref{ressit}. Let
  $p(x)\in\C[x]$ be a polynomial.
  There exists $i_0\in I$ depending on the matrix $A$ and on $p$ such that 
  \begin{equation*} 
    \tr_{\pi}(p(A)) = \tr_{i}(p(A_i))\qquad\forall
    i> i_0 .
  \end{equation*}
\end{lemma}
\begin{proof}
  Suppose $\pi$ is the inverse limit of the $\pi_i$.
  We follow \cite[2.6]{Lueck(1994c)}. Let
  $p(A)=\left(\sum_{g\in\pi}\lambda^{kl}_g g\right)_{k,l=1,\dots,d}$. Then
  \begin{equation*}
 \tr_{\pi}(p(A))=\sum_k\lambda^{kk}_1 \qquad\text{and}\qquad
  \tr_{\pi_i}(p(A_i))=
  \sum_k\sum_{g\in\ker p_i}\lambda_g^{kk} .
\end{equation*}
Since only finitely many
  of $\lambda_g^{ij}\ne
0$ and  $\pi$ is the inverse limit of the $\pi_i$  we find $i_0\in I$
such that $\lambda^{kk}_g\ne 0$ and
$g\in\ker p_{i_0}$ implies $g=1$. For $i> i_0$ the assertion is
true.

If $\pi$ is the direct limit we have chosen $A_{j_0}$ with
$p_{j_0}(A_{j_0})=A$. Then $p_{j_0}(p(A_{j_0}))=p(A)$. However, there
may be a $g\in\pi_{j_0}$ with nontrivial coefficient in $p(A_{j_0})$
with $p_{j_0}(g)=1$, and this means that the relevant traces may
differ. But still there are only finitely many $g$ with nontrivial
coefficient in $p(A_{j_0})$, and since $\pi$ is the direct limit of
$(\pi_i)_{i>j_0}$ we find $i_0$ such that $p_{i_0}(g)=1$ for $g\in
\pi_{i_0}$ with nontrivial coefficient in $p(A_{i_0})$ implies
$g=1$. Then the above reasoning shows
$\tr_{\pi_i}(p(A_i))=\tr_\pi(p(A))$ $\forall i>i_0$. Here for $i>j_0$
$A_i$ is the image of $A_{j_0}$ induced by the map $\pi_{j_0}\to \pi_i$.
\end{proof}

\section{Approximation properties and proofs of stability statements}
\label{sec:appr}

In this section, we use the information gathered so far to proof the
statements of the introduction, in particular Theorem
\ref{semiint_stab}. This is done by studying limits of
operators, therefore the same treatment yields approximation results
for $L^2$-Betti numbers. The precise statements and conditions are
given below.

We have already dealt with the passage to subgroups.
 
For the rest of this section, assume the following situation:
\begin{situation}\label{generalsit}
  The group $\pi$ is the direct or
  inverse limit of a directed system of groups $\pi_i$, or an amenable
  extension $U\to\pi\to\pi/U$ (write $\pi_i=U$ also in this case).

 As described in \ref{amsit} or \ref{ressit}, any
  matrix $\Delta$ over $\complexs\pi$ then gives rise to matrices $\Delta_i$ over $\pi_i$
  (after the choice of an inverse image in the case of a direct limit,
  and after the choice of an amenable exhaustion for amenable
  extensions). Without loss of generality we assume that $\Delta=A^*A$
  for another matrix $A$ over $\C\pi$ (this is explained in Remark
  \ref{oBdA}). We also get 
  spectral density functions
  $F_{A_i}(\lambda)$ defined using the group $\pi_i$ (remember that in
  the amenable case there is an additional normalisation).

  The problem now is to obtain information about $F_\Delta(\lambda)$
  from the family $F_{\Delta_i}(\lambda)$.

  In particular, we want to show that $F_{A_i}(0)$ converges to
  $F_A(0)$. (Translated to geometry this means that certain $L^2$-Betti
  numbers converge.)

  In short: we have
  \begin{itemize}
    \item A group $\pi$ and a matrix $\Delta\in M(d\times d,\complexs\pi)$
    \item A family $(\Delta_i){i\in I}$ of matrices over
      $\complexs\pi_i$ which approximate $\Delta$ ($I$ is a directed system)
    \item positive and normal trace functionals $\tr_i$ (on a von
      Neumann algebra which contains $\Delta_i$) which ares normalized in
      the following sense: If $\Delta=\id\in M(d\times
      d,\integers\pi)$ then $\tr_i(\Delta_i)=d$ $\forall i$.
    \item If $\Delta$ lives over $\integers\pi$, then $\Delta_i$ is a
      matrix over $\integers\pi_i$.
  \end{itemize}
\end{situation}

\begin{definition}
  Define
  \begin{equation*}
    \begin{split}
       \overline{F_\Delta}(\lambda) &:=
       \limsup_{i}F_{\Delta_i}(\lambda) ,\\
      \underline{F_\Delta}(\lambda) &:=
      \liminf_{i}F_{\Delta_m}(\lambda) .
    \end{split}
  \end{equation*}
Remember $\limsup_{i\in I}\{x_i\}=\inf_{i\in I} \{\sup_{j<i}\{x_j\}\}$.
\end{definition}

\begin{definition}
  Suppose $F:[0,\infty)\to \R$ is monotone increasing (e.g.\ a spectral
  density function). Then set
  \begin{equation*}
    F^+(\lambda):= \lim_{\epsilon\to 0^+} F(\lambda+\epsilon)
  \end{equation*}
  i.e.\ $F^+$ is the right continuous approximation of $F$. In
  particular, we have defined $\overline{F_\Delta}^+$ and
  $\underline{F_\Delta}^+$.
\end{definition}

\begin{remark}
  Note that
  by our definition a spectral densitiy function is right continuous,
  i.e.~unchanged if we perform this construction.
\end{remark}

To establish the first step in our program we have to establish the
following functional analytical lemma (compare \cite{Lueck(1994c)} or
\cite{Clair(1997)}):
\begin{lemma}\label{trIE}
  Let $\NeumannN$ be a finite von Neumann algebra with positive normal
  and normalized trace
  $\tr_\NeumannN$. Choose $\Delta\in M(d\times d,\NeumannN)$ positive
  and self-adjoint.\\
  If for a function $p_n:\R\to\R$
\begin{equation}\label{polIE}
 \chi_{[0,\lambda]}(x)\le p_n(x)\le 
  \frac{1}{n}\chi_{[0,K]}(x)+\chi_{[0,\lambda+1/n]}(x) \qquad\forall
  0\le x\le K
\end{equation}
 and if $\norm{\Delta}\le K$ then
 \begin{equation*}
   F_\Delta(\lambda)\le \tr_\NeumannN p_n(\Delta) \le
   \frac{1}{n}d + F_\Delta(\lambda+1/n) .
 \end{equation*}
  Here $\chi_{S}(x)$ is the characteristic function of the subset
  $S\subset \R$.
\end{lemma}
\begin{proof}
  This is a direct consequence of positivity of the trace, of the
  definition of spectral density functions and of the
  fact that $\tr_\NeumannN(1\in M(d\times d,\NeumannN))=d$ by the
  definition of a normalized trace.
\end{proof}

\begin{proposition}\label{appr}
  For every $\lambda\in\reals$ we have
  \begin{equation*} 
    \begin{split}
      &\overline{F_\Delta}(\lambda)\le F_\Delta(\lambda)=
      \underline{F_\Delta}^+(\lambda)\le \underline{F_\Delta}^+(\lambda) ,\\
      &
      F_\Delta(\lambda)
      =\underline{F_\Delta}^+(\lambda)=\overline{F_\Delta}^+(\lambda).
  \end{split}
  \end{equation*}
\end{proposition}
\begin{proof}
  The proof only depends on
     our key lemmas
    \ref{ambound}, \ref{amappr}, \ref{limbound}, \ref{limappr}. These say
    \begin{itemize}
    \item $\exists K\in\reals$ such that $\norm{\Delta_i}\le K$ $\forall
      i\in I$
    \item For every polynomial $p\in\complexs[x]$ we have
      $\tr_\pi(p(\Delta))=\lim_i\tr_i(p(\Delta_i))$.
    \end{itemize}

  For each $\lambda\in\reals$ choose polynomials $p_n\in\reals[x]$ such that
  inequality \eqref{polIE} is fulfilled. Note that by the first key
  lemma we find a uniform upper bound $K$ for the spectrum of all of
  the $\Delta_i$. Then by Lemma \ref{polIE}
  \begin{equation*}
    \begin{split}
      F_{\Delta_i}(\lambda) & \le \tr_i( p_n(\Delta_i)) \le
      F_{\Delta_i}(\lambda+\frac{1}{n}) + \frac{d}{n}
  \end{split}
  \end{equation*}
  We can take the limes inferior and superior
  and use the second key lemma to get
  \begin{equation*}
    \overline{F_\Delta}(\lambda) \le \tr_\pi(p_n(\Delta)) \le
    \underline{F_\Delta}(\lambda+\frac{1}{n}) +\frac{d}{n}.
  \end{equation*}
  Now we take the limes for $n\to\infty$. We use the fact that
  $\tr_\pi$ is normal and $p_n(\Delta)$ converges strongly inside a
  norm bounded set to $\chi_{[0,\lambda]}(\Delta)$. Therefore the
  convergence even is in the ultrastrong topology.

  This implies
  \begin{equation*}
    \overline{F_\Delta}(\lambda) \le F_\Delta(\lambda) \le
    \underline{F_\Delta}^+(\lambda) .
  \end{equation*}
  For $\epsilon>0$ we can now conclude since $\underline{F_\Delta}$
  and $\overline{F_\Delta}$ are monoton
\begin{equation*}
 F_\Delta(\lambda)\le \underline{F_\Delta}(\lambda+\epsilon)\le
 \overline{F_\Delta}(\lambda+\epsilon)\le F_\Delta(\lambda+\epsilon).
\end{equation*}
Taking the limit $\epsilon\to 0^+$ gives (since $F_\Delta$ is right
continuous)
\begin{equation*}
F_\Delta(\lambda)=\overline{F_\Delta}^+(\lambda)=
\underline{F_\Delta}^+(\lambda).
\end{equation*}
Therefore both of the inequalities are established.
\end{proof}

The next step is to proof convergence results without taking right
continuous approximations (at least for $\lambda=0$). We are able to
do this only under additional assumptions:
\begin{itemize}
\item From now on, $\Delta$ and therefore also $\Delta_i$ $\forall
  i\in I$ are matrices over the integral group ring.
\end{itemize}

The following statement is used
as start for the induction.
\begin{lemma}\label{trivgroup}
  The trivial group has semi-integral determinant.
\end{lemma}
\begin{proof}
  Take $\Delta\in M(d\times d,\integers)$ positive and
  self-adjoint. Then $\Det_1(\Delta)$ is the product of all nonzero
  eigenvalues and therefore the lowest nonzero coefficient in the
  characteristic polynomial. Therefore, it is an  integer $\ne 0$ and
  $\ln\Det_1(\Delta)\ge 0$.
\end{proof}

Now we give a proof of Theorem \ref{semiint_stab}  and prove the
corresponding approximation result.

\begin{theorem}\label{sintapr}
  Suppose 
  $\pi_i$  has
  ``seminitegral determinant'' $\forall i\in I$, then the same is true
  for $\pi$, and
  $\dim_\pi(\ker\Delta)=F_\Delta(0)= \lim_i F_{\Delta_i}(0)$.
\end{theorem}
\begin{proof}
 Choose $K\in\reals$ such that $K>\norm{\Delta}$ and $K>\norm{\Delta_i}$ $\forall
 i$. This is possible because of the key Lemma \ref{ambound} or
 \ref{limbound}. Then
  \begin{equation*}
    \ln\Det_{\pi_i}(\Delta_i)=\ln(K)(F_{\Delta_i}(K)-F_{\Delta_i}(0))
    -\int_{0^+}^{K}
    \frac{F_{\Delta_i}(\lambda)-F_{\Delta_i}(0)}{\lambda}\; d\lambda .
  \end{equation*}
  If this is (by assumption) $\ge 0$, then since
  $F_{\Delta_i}(K)=\tr_i(\1_d)=d$ by our normalisation
  \begin{equation*}
    \int_{0^+}^{K}\frac{F_{\Delta_i}(\lambda)-F_{\Delta_i}(0)}{\lambda}\;d\lambda \le
    \ln(K)(d-F_{\Delta_i}(0))\le \ln(K)d .
  \end{equation*}
 We want to establish the same estimate for
  $\Delta$. If $\epsilon>0$ then
  \begin{equation*}
    \begin{split}
       \int_{\epsilon}^K\frac{F_\Delta(\lambda)-F_\Delta(0)}{\lambda}\;d\lambda
      & = \int_{\epsilon}^K
      \frac{\underline{F_\Delta}^+(\lambda)-F_\Delta(0)}{\lambda}\;d\lambda
      =\int_{\epsilon}^K\frac{\underline{F_\Delta}(\lambda)-F_\Delta(0)}{\lambda}\\
      \intertext{(since  the integrand is bounded, the integral
        over the left continuous approximation is equal to the integral
        over the original function)}\\ 
      \le & \int_{\epsilon}^K
      \frac{\underline{F_\Delta}(\lambda)-\overline{F_\Delta}(0)}{\lambda}\\
      = & \int_\epsilon^K \frac{\liminf_i F_{\Delta_i}(\lambda)-\limsup_i
        F_{\Delta_i}(0)}{\lambda}\\
      \le &\int_\epsilon^K
      \frac{\liminf_i (F_{\Delta_i}(\lambda)-F_{\Delta_i}(0))}{\lambda}\\
      \le &\liminf_i\int_\epsilon^K
      \frac{F_{\Delta_i}(\lambda)-F_{\Delta_i}(0)}{\lambda} \le
      d\ln(K) .
  \end{split}
\end{equation*}
Since this holds for every $\epsilon>0$, we even have
\begin{equation*}
  \begin{split}
    \int_{0^+}^K \frac{F_\Delta(\lambda)-F_\Delta(0)}{\lambda} 
    \le & \int_{0^+}^K \frac{\underline{F_\Delta}(\lambda)-
      \overline{F_\Delta}(0)}{\lambda}\;d\lambda\\
    \le & \sup_{\epsilon>0}\liminf_{i} \int_{\epsilon}^K
    \frac{F_{\Delta_i}(\lambda)-F_{\Delta_i}(0)}{\lambda}\;d\lambda\\
    \le & d \ln(K) .
  \end{split}
\end{equation*}
The second integral would be infinite if $\lim_{\delta\to
  0}\underline{F_\Delta}(\delta)\ne \overline{F_\Delta(0)}$. It
follows that $\limsup_i
F_{\Delta_i}(0)=F_\Delta(0)$. Since we can play the same game for every
subnet of $I$, also $\liminf_i F_{\Delta_i}(0)=F_\Delta(0)$ i.e.\ the
approximation property is true.

For the estimate of the determinant note that in the above inequality 
\begin{equation*}
  \begin{split}
    & \sup_{\epsilon>0}\liminf_{m\to\infty}\int_{\epsilon}^K
    \frac{F_{\Delta_i}(\lambda)-F_{\Delta_i}(0)}{\lambda}\;d\lambda\\
    \le & \liminf_i \sup_{\epsilon>0}\int_\epsilon^K
    \frac{F_{\Delta_i}(\lambda)-F_{\Delta_i}(0)}{\lambda}\;d\lambda = \liminf_i \int_{0^+}^K
    \frac{F_{\Delta_i}(\lambda)-F_{\Delta_i}(0)}{\lambda}\;d\lambda\\
    \le & \ln(K)(d-F_{\Delta_i}(0)) .
  \end{split}
\end{equation*}
Therefore
\begin{equation*}
  \begin{split}
    \ln\Det_\pi(\Delta) &=\ln(K)(d-F_\Delta(0)) -\int_{0^+}^K
    \frac{F_\Delta(\lambda)-F_\Delta(0)}{\lambda}\;d\lambda\\
    &\ge \ln(K)(d-\lim_i F_{\Delta_i}(0)) - \liminf_i \int_{0^+}^K
    \frac{F_{\Delta_i}(\lambda)-F_{\Delta_i}(0)}{\lambda}\;d\lambda\\
    &= \limsup_i\left(\ln(K)(d-F_{\Delta_i}(0)) - \int_{0^+}^K
      \frac{F_{\Delta_i}(\lambda)-F_{\Delta_i}(0)}{\lambda}\;d\lambda\right)\\
    &= \limsup_i\ln\Det_{\pi_i}(\Delta_i) \ge 0. \qquad\qed
  \end{split}
\end{equation*}
\renewcommand{\qed}{}
\end{proof}

\begin{remark}
  The above reasoning does not show that the determinants converge. An
  unpublished example of L\"uck shows that this in general is wrong for
  matrices over the complex group ring of $\Z$. But he shows that the
  statement is true for matrices over $\Z[\Z]$. For other groups, the
  question is completely open.
\end{remark}

\begin{remark}
  The assumption ``semi-integral determinant'' is very
  strong. Originally, L\"uck as well as Dodziuk/Mathai and Clair use
  the following in some
  sense weaker
  property:
   a discrete group $\pi$ is \emph{spectrally sublogarithmic} if for $\Delta\in
    M(d\times d,\Z\pi)$ positive and self-adjoint
    \begin{equation*}
       F_\Delta(\lambda)-F_\Delta(0)\le d\frac{\ln R(\Delta)}{-\ln
         \lambda}\qquad \text{for $0<\lambda<1$}.
     \end{equation*}
   However, it is not clear how to find $R(\Delta)$ such that the
   property is stable under limits as well as under amenable
   extensions.
   
   As observed by Clair \cite{Clair(1997)}, instead of $1/\ln(x)$ any other
   function which is right continuous at zero would work equaly well.
\end{remark}

\begin{remark}
  We can establish the approximation results only under the
  assumption that the groups $\pi_i$ have good properties, e.g.~belong
  to the class $\RAgroups$. 
  It is interesting to note that for
  amenable groups every quotient
  group is amenable and  belongs to
  $\RAgroups$. It follows that Conjecture
  \ref{resappr} holds if $\pi_1(X)$ is amenable   without additional
  assumptions.

  If $\pi\in\RAgroups$ and $U\subgroup \pi$ then also
  $U\in\RAgroups$. Therefore, in case of an amenable exhaustion
  (i.e.~$\pi/U$ is amenable and we approximate using this fact) the
  conditions which imply convergence  are automatically fulfilled.
\end{remark}

\begin{remark}
  It would be possible to give geometrical interpretations of the more
  general approximation results. However, this seems to be very
  artificial and therefore is omitted here.
\end{remark}
\forget{
\subsection*{Geometrical interpretation of approximation results}
Since this geometric interpretation is not the focus of this paper and
seems to be of less importance, we only indicate briefly how this can
be done.

Amenable approximation (with groups $U\subgroup \pi$ and with the
changes indicated in Remark
\ref{geochange}) corresponds to an exhaustion of the finite
$\pi$-CW-complex $\tilde X$ by finite $U$-CW-complexes $\tilde X_n$ as
outlined in Dodziuk/Mathai \cite[Introduction]{Dodziuk-Mathai(1996a)}. Residual
approximation is explained in the introduction. The  limit
cases corresponds to sequences
of Hilbert bundles (fiber $l^2(\pi_i)$) over a fixed space $X$,
converging to the Hilbert bundle with
fiber $l^2(\pi)$.
}
\subsection{Proof of Theorem \ref{hominvstab}}

It remains to show that Whitehead trivial
determinant is stable under direct and inverse limits.
\begin{proof}[Proof of Theorem \ref{hominvstab}]
  We are still in the situation
  described at the beginning of this section, and, in addition, we assume that
  $\Delta$ is invertible in 
  $M(d\times d,\Z \pi)$ with inverse $B\in M(d\times d,\Z\pi)$. In
  case $\pi$ is an inverse limit $\Delta_i$ and $B_i$ are images of
  projections of $\Delta$ and $B$, and therefore remain inverse to
  each other.

  In case $\pi$ is a direct limit, we first lifted $\Delta$ to some
  $\Delta_{j_0}$. We may assume that we also can lift $B$ to
  $B_{j_0}$. Then $\Delta_{j_0}B_{j_0}$ is mapped to the identity over
  $\pi$. Since it has only finitely many nonzero coefficients, there
  is $j_1$ such that the image of $\Delta_{j_0}B_{j_0}$ over $\pi_{j_1}$
  already is the identity, and similar for
  $B_{j_0}\Delta_{j_0}$. Therefore, we may assume that the lifts
  $\Delta_{j_0}$ and $B_{j_0}$ are inverse to each other. The same is
  then true for $\Delta_i$ and $B_i$ for $i>j_0$, i.e.~$\Delta_i$
  represents an element in $Wh(\pi_i)$.

  By assumption
  $\ln\Det_{\pi_i}(\Delta_i^*\Delta_i)=0$. Note that the proof of
  Theorem \ref{sintapr} applies to our situation and we conclude
  $\ln\Det_\pi(\Delta^*\Delta)\ge 0$.
  Since $\Delta\in Wh(\pi)$ was arbitrary, Theorem \ref{T20} implies
  the result.
\end{proof}

\begin{remark}
  It is not possible to proceed along similar lines in the case that
  $U$ has an amenable quotient (even if $\pi$ is amenable
  itself). The problem is that we approximate the matrix $\Delta$
  (over $\Z\pi$) by matrices over $\Z U$ of larger and larger
  dimension. One can show that these matrices are invertible over
  $\NeumannN U$, if $\Delta$ itself was invertible. However, even if
  the inverse of $\Delta$ is a matrix over $\Z\pi$, in general this is
  not true for the approximating matrices over $\Z U$.
\end{remark}

This finishes the proof of Theorems \ref{hominvstab},
\ref{semiint_stab} and \ref{main}. We proceed with some side remarks.

\section{Complex approximation}
In this section, we will adress the question wether the approximation
results we have obtained in section \ref{sec:appr} are valid not only
for matrices over the
integral group ring, but also over the complex group ring. In
particular, we adopt Situation \ref{generalsit}: A group $\pi$ is
approximated by groups $\pi_i$, and a matrix $\Delta=A^*A$ by matrices
$\Delta_i$.

Also, we try to relate the approximation problem to the Atiyah
conjecture, which says
\begin{conjecture}
  Suppose $\pi$ is torsion free. Then $\dim_\pi(\ker A)\in\integers$
  whenever $A\in M(d\times d,\complexs\pi)$. 
\end{conjecture}
This conjecture is true for abelian groups, free
groups, and for extensions with amenable quotient (compare Linnell
\cite{Linnel(1993)})

Essentially, we will give a positive answer to our question only for
free abelian groups.
We start with a general observation.
\begin{lemma}\label{0appr}
  Suppose in the situation \ref{generalsit} that 
  $\ker \Delta=0$. Then the
  approximation result holds without integrality assumptions: $\lim_i F_{\Delta_i}(0)\, = 0 = \dim_\pi(\ker
  \Delta)$.

  More generally, if $\lambda$ is not an eigenvalue of $\Delta$ then
  $\lim_i F_{\Delta_i}(\lambda)=F_\Delta(\lambda)$.
\end{lemma}
\begin{proof}
  We know that $\underline{F}_\Delta^+(x)=F_\Delta(x)$ for every $x\in\R$. If
  $\underline{F}_\Delta(\lambda)\le F_\Delta(\lambda)-\epsilon$ then
  $\underline{F}_\Delta^+(x)\le F_\Delta(\lambda)-\epsilon$ for every
  $x<\lambda$, i.e.~$F_\Delta(x)\le F_\Delta(\lambda)-\epsilon$ $\forall
  x<\lambda$. By assumption, the eigenspace of $\Delta$ to $\lambda$ is
  trivial, therefore $F_\Delta$ is continuous at $\lambda$ and $\epsilon$
  can only be zero.
\forget{  (For the last continuity statement::
  $F_\Delta(\lambda)-F_\Delta(\lambda-\delta)=
  \tr_\pi\chi_{[\lambda-\delta,\lambda)}(\Delta)$. The family of characteristic
  functions converges pointwise to the zero function, therefore the
  corresponding spectral projections converge strongly to zero. Since
  their norm is bounded by $1$, the convergence even is ultrastrongly
  and the regularized traces also converge to zero.) 
}
\end{proof}

\begin{proposition}\label{1dimappr}
  If $\pi$ is torsion free and  fulfills the Atiyah conjecture for
  $1\times 1$-matrices, then the above approximation
  result $\lim_i F_{\Delta_i}(0)=F_\Delta(0)$ holds
  $\forall\Delta=A^*A\in \complexs\pi$.
\end{proposition}
\begin{proof}
  Take $0\ne \Delta\in\complexs\pi$. By assumption, since $\ker \Delta\ne
  l^2(\pi)$, $\ker
  \Delta=0$. Conclude with lemma \ref{0appr}.
\end{proof} 

\begin{proposition}\label{comp_appr}
  If $\pi$ is free abelian  and a (subgroup of an) inverse limit
  ---e.g.~if we approximate $\pi$ residually---
 then the approximation
  result holds $\forall\Delta\in\complexs\pi$.
\end{proposition}
\begin{proof}
   Since all properties are stable under directed unions, it suffices
   to consider a finitely generated free abelian group
   $\pi\cong\Z^n$. Proposition \ref{1dimappr} shows that the statement
   holds for $1\times 1$-matrices. We will use commutativity to reduce
   the general case to the one-dimension case. Embed $\C\pi$ into its
   ring of fractions. Let $A\in M(d\times d,\C[\Z^n])$. Linear algebra
   tells us that we find $X,Y\in Gl(d\times d,\C(\Z^n))$ such that
   $A=X\diag(1,\dots,1,0,\dots,0)Y$. Collecting the denominators in
   $X$ and $Y$ we find $0\ne c\in \C[\Z^n]$ and $X',Y'\in M(d\times
   d,\C[\Z^n])\cap Gl(d\times d,\C\pi)$ such that 
   \begin{equation*}
     cA = X'\diag(1,\dots,1,0,\dots,0)Y'
   \end{equation*}
   The corresponding equation holds after passage to
   matrices over $\pi_i$. Since $c$, $X'$ and $Y'$ have trivial kernel, lemma
   \ref{0appr} implies the desired convergence result.
\end{proof}

\begin{remark}
  With a little bit more effort, one can get similar approximation
  results also in the other two contexts we are studying, in
  particular for amenable exhaustions of $\Z^n$.
\end{remark}

\begin{remark}
  Approximation with complex coefficients implies that not only the dimensions
  or the kernels but of the eigenspaces to every complex number
  converge, since we simply replace $A$ by $A-\lambda$. Together with
  the fact that the right continuous $limsup$ of the $F_{\Delta_i}$ is
  the spectral density function of $\Delta$, we
  have convergence of the spectral density function at every point.
\end{remark}

\section{Quotients}
\label{sec:quot}

To enlarge the class $\RAgroups$, it is important to find other
operations under which our main properties: determinant class and
semi-integrality are inherited.

We indicate just one partial result:
\begin{proposition}
  Suppose $1\to F\to \pi\to Q\stackrel{p}{\to} 1$ is an extension of groups and
  $\abs{F}<\infty$, and $\pi$ is of determinant class. Then also $Q$
  is of determinant class.
\end{proposition}
\begin{proof}
  We only indicate the proof, which was discussed with M.~Farber
  during a conference in
  Oberwolfach, and which uses the theory of virtual characters
  of  Farber \cite{Farber(1997)}. $l^2(\pi)$
  corresponds to the Dirac character $\delta_1$. The
  representations $V_k$ of the finite group $F$ give rise to
  characters $\chi_k$ of $\pi$ with support contained in $F$. Since
  $l^2(F)=\bigoplus \mu_k V_k$, $\delta_1=\frac{1}{\abs{F}}\sum
  \mu_k\chi_k$ with $\mu_k>0$. Since the operator $\Delta$ we are interested
  in arises from $c^*\tensor \id$ on $C^*\tensor l^2(\pi)$,
  $F_\Delta(\lambda)= \sum
  \frac{\mu_k}{\abs{F}} F^{\chi_k}_{\Delta}(\lambda)$ (compare
  \cite[7.2]{Farber(1997)}). Now the trivial
  representation $V_1$ of $F$ corresponds to the quotient 
  representation $l^2(Q)$, and $F^{\chi_1}_{\Delta}$ is just
  $F^Q_{p(\Delta)}$. By assumption $\int_{0+}^\infty
  \ln(\lambda)\; dF_\Delta(\lambda)>-\infty$. Since $\int\ln(\lambda)\;
  dF^{\chi_k}_{\Delta}(\lambda)<\infty$ $\forall k$, it follows in particular
 \begin{equation*}
    \int_{0^+}^\infty \ln(\lambda)\; dF^Q_{p(\Delta)}(\lambda)> -\infty.
  \end{equation*}
  Since $p$ is surjective, this is true for every matrix over
  $\integers\pi$ we have to consider. This concludes the proof.
\end{proof}

\begin{remark}
  If we have an extension $1\to\integers^n\to\pi\stackrel{p}{\to}Q\to
  1$, we can not write the character of $l^2(\pi)$ as a direct sum,
  but as a direct integral (over the dual space
  $\widehat{\integers^n}=T^n$). Then
  $F_\Delta(\lambda)=\int_{T^n}F^{\chi_\eta}_{\Delta}(\lambda)\;d\eta$, and
  $F^{\chi_1}_{\Delta}(\lambda)=F^Q_{p(\Delta)}(\lambda)$. After
  establishing an appropriate continuity property, we can conclude as
  above that  if $\pi$ is of
  determinant class then
\begin{equation*} \int_{0^+}^\infty
  \ln(\lambda)\;dF^Q_{p(\Delta)}(\lambda)>-\infty,
\end{equation*}
i.e.~also $Q$ is of determinant class.
  
\end{remark}

\forget{
The method does not apply to conclude results about semi-integral
determinant. On the other hand, the above splittings suggest that one
should study other representation besides the regular one, even if one
is only interested in results concerning the latter. This was startet
by Farber \cite{Farber(1997)}, but our more general approach should
give more information also here.
}
\ifx\undefined\bysame
\newcommand{\bysame}{\leavevmode\hbox to3em{\hrulefill}\,}
\fi

\end{document}